\documentclass[12pt]{amsart}

\usepackage{amsmath}
\usepackage{amssymb,amscd,pb-diagram,mathrsfs}
\usepackage{fullpage}
\usepackage[all]{xy}

\newcount\timehh\newcount\timemm
\timehh=\time
\divide\timehh by 60 \timemm=\time
\def\C{{ \mathbb C }}

\def\Z{{   \mathbb Z }}
\def\Q{{\mathbb Q}}

\def\F{{\mathbb F}}
\def\1{{\bf 1}}
\def\G{Sl_2(\Z)}

\def\oo{{\mathcal O}}

\def\ee{{\mathscr F}}
\def\kk{{\mathscr K}}

\def\ph{{\varphi}}
\def\Ph{{\Phi}}

\def\cP{{\mathcal P}}

\def\Aut{{\rm Aut~~}}
\def\compos{{\cup}}

\def\sig{ \sigma }

\def\ord{{\rm ord}}

\def\G{{\mathcal G}}

\def\m{{\frak m}}
\def\P{{\mathbb P}}
\def\p{{\frak p}}

\def\rd{{ \mathop{{\rm {rd}}} }}
\def\Gal{{ \mathop{{\rm {Gal}}} }}

\def\disc{{ \mathop{{\rm {disc}}} }}

\def\mod{\text{mod}}

\def\Gl{{\mathrm{GL}}}

\def\Res{{\rm Res}}

\def\c2{\chi_2}

\def\hfl#1#2{\smash{\mathop{\hbox to
12mm{\rightarrowfill}}\limits_{\scriptstyle #1}^{\scriptstyle #2}}}

\def\C{{ \mathbb C }}

\def\Z{{   \mathbb Z }}
\def\Q{{\mathbb Q}}

\def\F{{\mathbb F}}
\def\G{Sl_2(\Z)}

\def\o{{\mathcal O}}

\def\fn#1{\ph^{\circ #1}}

\def\G{{G}}

\def\m{{\mathfrak m}}
\def\P{{\mathbb P}}
\def\p{{\mathfrak p}}

\def\rd{{ \mathop{{\rm {rd}}} }}
\def\Gal{\text{Gal}}

\def\disc{{ \mathop{{\rm {disc}}} }}

\def\deg{{ \mathop{{\rm {deg}}} }}

\def\mod{\text{mod}}

\def\RR{{\mathcal R}}
\def\B{{\mathcal B}}
\def\PP{{\mathbb P}}
\def\M{{\mathscr M}}
\def\diff{{\mathfrak d}}
\def\PC{{\mathcal P}}

\def\m{    { {{\mathfrak t}}_0 }              }

\newtheorem{theorem}{Theorem}[section]
\newtheorem{lemma}[theorem]{Lemma}
\newtheorem{corollary}[theorem]{Corollary}
\newtheorem{proposition}[theorem]{Proposition}

\theoremstyle{definition}
\newtheorem{definition}[theorem]{Definition}
\newtheorem{example}[theorem]{Example}

\newtheorem{question}[theorem]{Question}
\theoremstyle{remark}

\newtheorem*{remark}{Remark}

\newenvironment{Proof}{\removelastskip\par\medskip
\noindent{\em Proof.} \rm}{\penalty-20\null\hfill$\square$\par\medbreak}

\begin{document}

\title[Finitely Ramified Iterated Extensions]
{Finitely Ramified Iterated Extensions}

\author{Wayne Aitken}
\email{waitken@csusm.edu} 
\address{Department of Mathematics, CSUSM,
San Marcos CA 92096 USA}

\author{Farshid Hajir}
\email{hajir@math.umass.edu}
\address{Department of Mathematics \& Statistics, 
University of Massachusetts,
	Amherst, MA 01003-9318 USA}

\author{Christian Maire}
\email{ christian.maire@univ-tlse2.fr}
\address{GRIMM,
U.F.R.~S.E.S.,
Universit\'e Toulouse 2,
5, all\'ees Antonio Machado,
31058 Toulouse c\'edex
France   }



\thanks{Hajir's work was partially supported by the National Science
Foundation under Grant No. DMS-0226869.}

\begin{abstract}
Let $K$ be a number field, $t$ a parameter, $F=K(t)$, and $\ph(x)\in
K[x]$ a polynomial of degree $d\geq 2$.  The polynomial $\Ph_n(x,t) =
\fn{n}(x)-t \in F[x]$, where $\fn{n}=\ph\circ \ph \circ \cdots \circ
\ph$ is the $n$-fold iterate of $\ph$, is absolutely irreducible over
$F$; we compute a recursion for its discriminant.  Let $\ee_\ph$ be
the field obtained by adjoining to $F$ all roots (in a fixed
$\overline{F}$) of $\Ph_n(x,t)$ for all $n\geq 1$; its Galois group
$\Gal(\ee_\ph/F)$ is the iterated monodromy group of $\ph$.  The
iterated extension $\ee_\ph$ is finitely ramified over $F$ if and only
if $\ph$ is post-critically finite.  We show that, moreover, for
post-critically finite $\ph$, every specialization of $\ee_\ph/F$ at
$t=t_0 \in K$ is finitely ramified over $K$, pointing to the
possibility of studying Galois groups of number fields with restricted
ramification via tree representations associated to iterated monodromy
groups of post-critically finite polynomials.  We discuss the wildness
of ramification in some of these representations, describe prime
decomposition in terms of certain finite graphs, and also give some
examples of monog\`ene number fields that arise from the construction.
\\ \begin{center} August 12, 2004\end{center}
\end{abstract}

\maketitle

\section{Introduction}\label{sec:intro}

Let $p$ be a prime number, $K$ a number field, and $S$ a finite set of
places of $K$.  Let $K_S$ be the compositum of all extensions of $K$
(in a fixed algebraic closure $\overline{K}$) which are unramified
outside $S$, and put $G_{K,S}=\Gal(K_S/K)$ for its Galois group.
These {\em arithmetic fundamental groups} play a very important role
in number theory.  Algebraic geometry provides the most fruitful known
source of information concerning these groups.  Namely, given a smooth
projective variety $X/K$, the $p$-adic \'etale cohomology groups of
$X$ are finite-dimensional vector spaces over $\Q_p$ equipped with an
action of $G_{K,S}$ where $S$ consists of the primes of bad reduction
for $X/K$ together with the primes of $K$ of residue characteristic
$p$.  The richness of this action can be judged, for example, by the
intimate relationships between algebraic geometry and the theory of
automorphic forms which it mediates.


For this and many other reasons, it would be difficult to overstate
the importance of
these $p$-adic Galois representations.
Nonetheless, linear $p$-adic groups simply form too restrictive a
class of groups to capture {\em all} Galois-theoretic information, and
some important conjectures in the subject, notably the Fontaine-Mazur
conjecture \cite{fm} (to mention only one, see the discussion in
section \ref{sec:conjectures}), point specifically toward the kind of
information inside arithmetic fundamental groups which {cannot} be
captured by finite-dimensional $p$-adic representations.

In this work, we discuss a method for studying finitely ramified
extensions of number fields via {\em arithmetic dynamical systems on}
${\PP^1}$.  At least conjecturally, this method provides a vista on a
part of $G_{K,S}$ invisible to $p$-adic representations.
We now sketch the construction, which is quite elementary.  Let $K$ be
a perfect field, and suppose $\ph\in K[x]$ is a polynomial of degree
$d\geq 1$ such that its derivative $\ph'$ is not identically $0$ in
$K[x]$.  For each $n\geq 0$, let $\fn{n}$ be the $n$-fold iterate of
$\ph$, i.e.  $\fn0(x)=x$ and
$\fn{n+1}(x)=\ph(\fn{n}(x))=\fn{n}(\ph(x))$ for $n\geq 0$.  Let $t$ be
a parameter for $\PP^1_{/K}$ with function field $F=K(t)$.  We are
interested in the tower of branched covers of $\PP^1$ given by
\begin{equation}\label{pn}
\Ph_n(x,t) = \fn{n}(x) - t \in F[x],
\end{equation}
as well as extensions of $K$ obtained by adjoining roots of its
specializations at arbitrary $t_0 \in K$.  The variable-separated
polynomial $\Ph_n(x,t)$ is clearly absolutely irreducible (since it is
linear in $t$) and of degree $d^n$ in $F[x]$; it is separable over $F$
by the assumption that $\ph'$ is not identically 0.

Fix an algebraic closure $\overline F$ of $F$, and let $\overline{K}$
be the algebraic closure of $K$ determined by this choice, i.e. the
subfield of $\overline{F}$ consisting of elements algebraic over $K$.
For $n\geq 0$, let $T_{\ph,n}$ be the set of roots in $\overline{F}$
of $\Ph_n(x,t)$; it has cardinality $d^n$.  We denote by
$T_\ph$  the $d$-regular rooted tree whose vertex set is 
$\cup_{n\geq 0} T_{\ph,n}$, 
and whose edges connect $v$ to $w$ exactly when $\ph(v)=w$; its root
(at ground level) is $t$.

Let us choose and fix an {\em end} $\xi=(\xi_0,\xi_1, \xi_2,
\cdots)$ of this tree; in other words, we have $\ph(\xi_1)=\xi_0=t$
and $\ph(\xi_{n+1})=\xi_{n}$ for $n\geq 1$.  For each $n\geq 1$, we
consider the field $F_n = F(\xi_n)\simeq F[x]/(\Ph_n)$ and its Galois
closure $\ee_n=F(T_{\ph,n})$ over $F$.  Let $\oo_{\ee_n}$ be the 
integral closure of $K[t]$ in $\ee_n$.
Corresponding to each $t_0 \in
K$, we may fix compatible {\em specialization maps}
$\sig_{n,t_0}:\oo_{\ee_n} \rightarrow \overline{K}$ with image $\kk_{n,t_0}$,
a normal extension field of $K$ and put $\xi_n|_{t_0}=\sig_{n,t_0}(\xi_n)$
for the correspoding compatible system of roots of $\Ph_n(x,t_0)$.
We denote by $K_{n,t_0}$ the image of the restriction of $\sig_{n,t_0}$ to
$\oo_{F_n}$.
We refer the reader to subsection \ref{subsec:specialization} for more
details, but we should emphasize here that $\Ph_n(x,t_0)$ is not
necessarily irreducible over $K$; hence, although $\kk_{n,t_0}$
depends only on $\ph,n$ and $t_0$, the isomorphism class of 
$K_{n,t_0}$ depends {\em a priori}
on the choice of $\xi$ as well as on the choice of compatible
$\sig_{n,t_0}$.  Also, the Galois closure of $K_{n,t_0}/K$ is contained in,
but possibly distinct from, $\kk_{n,t_0}$.  Nonetheless, unless stated
otherwise, $\xi$ and $t_0$ are arbitrary but fixed, and in this case
we will usually not decorate
$K_{n,t_0}$ with $\xi$ and occasionally we may write simply $\kk_n,
K_n$ instead of $\kk_{n,t_0}, K_{n,t_0}$.

Taking the compositum over all $n\geq 1$, we obtain the {\em iterated
extension} $F_\ph=\compos_n F_n$ attached to $\ph$, with Galois
closure $\ee_\ph=\compos_n \ee_n$ over $F$.  Similarly for each $t_0
\in K$, we obtain a specialized iterated extension
$K_{\ph,t_0}=\compos_n K_{n,t_0}$ with Galois closure over $K$
contained in $\kk_{\ph,t_0}=\compos_n \kk_{n,t_0}$.
We put $\M_\ph = \Gal(\ee_\ph/F)$ for the {\em iterated monodromy
group} of $\ph$ and for $t_0 \in K$, we denote by
$\M_{\ph,t_0}=\Gal(\kk_{\ph,t_0}/K)$ its specialization at $t_0$.  The
group $\M_\ph$ has a natural and faithful action on the tree $T_\ph$,
hence comes equipped with a rooted tree representation $\M_\ph
\hookrightarrow \Aut T_\ph$.  For more on rooted trees and  iterated monodromy
groups (in a more general context, in fact), see Nekrashevych
\cite{nar} as well as Bartholdi-Grigorchuk-Nekrashevych \cite{bgn}.

Since we are interested in finitely ramified towers (meaning those
where only finitely many places of the base field are ramified), we
need to answer the following question: Which polynomials $\ph$ have
the property that the corresponding iterated tower $F_\ph/F$, as well
as all of its specialziations $K_{\ph,t_0}/K$, are finitely ramified?

We first recall some standard terminology from
polynomial dynamics.  We say that $z\in \overline{F}$ is {\em
periodic} for $\ph$ if $\fn{n}(z) = z$ for some $n\geq 1$.  Moreover,
$y\in \overline{F}$ is {\em preperiodic} for $\ph$ if for some $m\geq
0$, $\fn{m}(y)$ is periodic for $\ph$; equivalently, $y$ is
preperiodic for $\ph$ means that $\{ \fn{n}(y):n\geq 0\}$, i.e. the
orbit of $y$ under the iterates of $\ph$, is a finite set.
We put
$$
\RR_\ph := \{ r\in \overline{K}:\ph'(r)=0\}, \qquad \B_\ph := \{ \ph(r) : r\in
\RR_\ph\}
$$ for the set of affine {\em ramification} and {\em branch} points,
respectively.  The elements of $\RR_\ph$ and $\B_\ph$ are also the
{\em critical points}, respectively {\em critical values} of $\ph$.
The polynomial $\ph$ is called {\em post-critically finite} if
every member of $\RR_\ph$ is a preperiodic point for $\ph$.  In other
words, $\ph$ is post-critically finite exactly when the {\em
post-critical set} $\PC_\ph$, i.e. the union of the orbits of critical
points under the iterates of $\ph$, is a finite set.  It has long been
known that the post-critical set plays a crucial role in the dynamics
of the polynomial.  Indeed, the class of dynamical systems
corresponding to post-critically finite polynomials is a well-studied
one, having gained prominence following a celebrated theorem of
Thurston; see, for example, Douady-Hubbard \cite{dh},
Bielefed-Fisher-Hubbard \cite{bfh}, as well as the papers by Poirier
\cite{poirier}, Pilgrim \cite{pilgrim}, and Pakovich \cite{pakovich};
the latter two concern the connection with actions of
$\Gal(\overline{\Q}/\Q)$ on certain finite trees.

Our starting point is the following characterization of finitely
ramified iterated extensions.

\begin{theorem}\label{thm:1}
The iterated tower of function fields
$\ee_\ph/F$ is finitely ramified if and only if $\ph$
is post-critically finite.  If $\ph$ is post-critically finite, 
every specialization $\kk_{\ph,t_0}/K$ of this tower is 
finitely ramified.  
\end{theorem}
The first assertion of the
theorem is clear geometrically since $\B_{\fn{n}}= \B_{\ph}\cup
\ph(\B_\ph) \cup \cdots \cup \fn{n-1}(\B_{\ph})$.  The second
assertion, however, is not a formal consequence of the first, since
any element of $K$ is a unit in $F$; the proof we give for it proceeds
via Proposition \ref{prop:recur}, where we calculate a recurrence for
the discriminant of $\Ph_n(x,t)$ (valid for an arbitrary polynomial
$\ph$), giving a more precise version of the theorem.  The proof of 
Proposition \ref{prop:recur} rests on a Riemann-Hurwitz genus formula
for polynomials \cite{simon}.  

Now let us return to the case of a number field $K$.  For each
post-critically finite $\ph\in K[x]$, and each $t_0 \in K\setminus
\cP_\ph$, Theorem \ref{thm:1} provides a surjection
$\rho_{\ph,t_0}:\G_{K,S} \twoheadrightarrow \M_{\ph,t_0}$ for an
effectively determined finite set $S=S_{\ph,t_0}$ of places of $K$
(see Definition \ref{def:S} and Corollary \ref{cor:imr}).  We call
$\rho_{\ph,t_0}$ the {\em iterated monodromy representation} attached
to $\ph$ and $t_0$.


The study of automorphism groups of rooted trees is a relatively new
and quite active topic in group theory (see \cite{bort}, \cite{nar},
and \cite{bgn}).  The structure of non-abelian subgroups of these
automorphism groups appears to be quite different from that of linear
$p$-adic groups (see the papers just cited as well as Bux-Perez
\cite{bux-perez}).  The natural action of iterated monodromy groups on
rooted trees leads us to the expectation that 
iterated monodromy representations $\rho_{\ph,t_0}$
attached to post-critically finite polynomials $\ph \in K[x]$ have the
potential of revealing aspects of arithmetic fundamental groups which
are not visible to $p$-adic representations; see the discussion in
section \ref{sec:conjectures} as well as Boston's preprint
\cite{boston-tr}, where tree representations are suggested as the
proper framework for studying finitely ramified tame extensions.

Since all finitely ramified $p$-adic Galois representations with
infinite image are expected, by a conjecture of Fontaine and Mazur, to
be wildly (even deeply) ramified at some prime of residue
characteristic $p$, an immediate question is what can be said about
the presence of wild ramification in specialized iterated extensions
$\kk_{\ph,t_0}/K$.  Experimentation leads to the expectation that
generically the primes of residue characteristic dividing $d$ ramify
deeply in $\kk_{\ph,t_0}/K$.  For example, if $\ph(x)=x^d$ with $d>1$
and $K=\Q$, then for all $t_0 \in \Q$, the extensions
$\kk_{\ph,t_0}/K$ are deeply ramified at all $p$ dividing $d$.  (See,
however, Questions \ref{ques:1} and \ref{ques:2} in \S
\ref{sec:conjectures}).

Under an assumption of good reduction for $\ph$, we prove a partial
result toward this expectation, namely for integral $t_0$, we estimate
from below the power of $p$ dividing the discriminant of
$\Ph_n(x,t_0)$.
To be precise, in \S \ref{sec:wild}, we will prove
the following theorem.
\begin{theorem}\label{thm:wild1}
Let $K$ be a number field.  Suppose $\ph\in K[x]$ is post-critically
finite, has degree divisible by $p$, and has good reduction at a
valuation $v$ of residue characteristic $p$, i.e. $\ph$ has $v$-integral
coefficients with $v$-unital leading coefficient.
Then for any $t_0 \in \oo_K$,
$$v ( \disc~ \Ph_n(x,t_0) ) \geq n d^n v(p).$$
\end{theorem}  
Assuming $\Ph_n(x,t_0)$ is $K$-irreducible for all $n$, this estimate
shows that the tower of rings $\oo_K[\xi_n|_{t_0}]$, where $(\xi_n|_{t_0})$
is a compatible sequence of roots of $\Ph_n(x,t_0)$, is wildly
ramified at $p$.  Note that
$\oo_K[\xi_n|_{t_0}]$ is
an {\em order} inside the maximal order of $K(\xi_n|_{t_0})$; it is the
discriminant of the latter which is our primary interest, but the theorem
estimates the discriminant of the former.  This is one sense in which
the above theorem is only a partial answer to our
question about wild ramification in iterated extensions.  On the other
hand, in section \ref{sec:quad},
we illustrate with the tower corresponding to $\ph(x)=x^2-2$, the
possibility that the orders $\Z[\xi_n|_{t_0}]$ (for a large set of $t_0
\in \Z$) are maximal, giving examples of monog\`ene number fields.

The organization of this article is as follows.  In \S
\ref{sec:prelim}, we outline some preliminary facts regarding
post-critically finite polynomials, including a classification of the
very simplest examples for each degree, namely those that are {\em
critically fixed} (every critical point is fixed, also known as {\em
conservative}) and {\em simply ramified} (every non-trivial
ramification index is $2$).  In \S\S \ref{sec:pcf} and \ref{sec:wild},
we prove Theorems \ref{thm:1} \ref{thm:wild1}), respectively. In \S
\ref{sec:decomposition}, we describe the decomposition of unramified
primes in iterated towers in terms of simple properties of certain
finite graphs.  In \S \ref{sec:quad}, we study the quadratic case in
more detail, obtaining a recursion for writing down post-critically
finite quadratic polynomials, which give number fields of independent
interest; we also discuss the example $x^2-2$ in detail, proving
monogenicity of certain number fields.  Finally, in \S
\ref{sec:conjectures}, we outline a number of questions and open
problems.

\

\noindent{\bf Acknowledgments.}  
Work on this project began in June 2003 during visits by WA and CM to
Amherst supported partially by NSF grant DMS-0226869.
A few months later, the preprint \cite{nar} of V.~Nekrashevych was
posted on \verb?arxiv.org?; in it, iterated monodromy groups from a
more general perspective are studied.  We learned from \cite{nar} that
the construction (\ref{pn}) we had been studying had earlier been
suggested by R.~Pink (private communication).  A special case also
occurs in Boston's preprint \cite{boston-tr}.  CM thanks R.~Pink for a
helpful conversation.  CM and FH would like to thank
E.~Bayer-Fluckiger for financial support of their visits to EPFL,
Lausanne in 2004 where a portion of this work was completed.  We are
grateful to R.~Benedetto and P.~Gunnells for useful discussions.

\section{Preliminaries}\label{sec:prelim}

\subsection{The branched cover $\fn{n}:\PP^1 \rightarrow \PP^1$}

In this section, $K$ is a perfect field and $\ph(x)=a_dx^d + \ldots +
a_0 \in K[x]$ is a polynomial of degree $d\geq 1$ whose derivative
$\ph'$ is not identically $0$.
We maintain all other notation introduced in \S 1.


Thinking of $\ph$ as a branched cover $\PP^1\rightarrow \PP^1$ of
degree $d$, the singular fibers are those of cardinality less than
$d$.  Leaving aside $\infty$ which is totally ramified, the points in
a singular fiber (the ramification points) are exactly the critical
points, i.e. the roots of $\ph'$: writing
$\ph(x)-\ph(r)=(x-r)\psi_r(x)$ for any $r\in \overline{K}$, we have
$\ph'(r)=\psi_r(r)$ hence $(x-r)^2$ divides $\ph(x)-\ph(r)$ if and
only if $\ph'(r)=0$.  The critical values (the images under $\ph$ of
the critical points), are the points having a singular fiber, i.e. they are
exactly the branch points.
In algebraic language, $\beta\in \overline{K}$ is in $\B_\ph$ if and
only if $\ph(y)-\beta$ has a multiple root, which happens if and only if
$\disc_y(\ph(y)-x)$ has $x=\beta$ as a root. In other words, $\beta$ is
a branch point if and only if the system $$\ph(y)=\beta, \qquad
\ph'(y)=0$$ has a common root $y=r$, and these roots are the
ramification points above $\beta$.  We could adopt the
convention that $\RR_\ph$ and $\B_\ph$ are
``multisets'' where each critical point or critical value occurs
according to the multiplicity of the corresponding roots of $\ph'$, but
to avoid confusion, we will be explicit about the multiplicities by
writing 
\begin{equation}\label{mr}
\ph'(x)= da_d \prod_{r \in \RR_\ph} (x-r)^{m_r}
\end{equation}
 and putting, for
$\beta \in \B_\ph$,
\begin{equation}\label{Mbeta}
M_\beta = \sum_{r\in \RR_\ph, \ph(r)=\beta} m_r.
\end{equation}

\begin{lemma}\label{lem:integral}
For each $n\geq 1$, $\Ph_n(x,t)$ is separable and absolutely 
irreducible over $F$.  The ring $K[\xi_n,t]$ is integrally closed
(in its fraction field $F_n$).
\end{lemma}
\begin{proof}
All of this follows essentially from the fact that $\partial_t
\Ph_n(x,t)=1$ never vanishes. The reader can easily check the absolute
irreducibility of $\Ph_n$. For separability, assume that $\Ph_n(x,t)$
has a multiple root, $\xi_n$ say.  Then $\xi_n$ is a root of $\partial_x
\Ph_n(x,t)=(\fn{n})'(x).$ Since $\ph'$ is not identically $0$, neither
is $(\fn{n})'$, and so $\xi_n$ is algebraic over $K$, and then so is
$t=\fn{n}(\xi_n)$, a contradiction.  Note that if $\ph'\equiv 0$, then
$\Ph_n(x,t)$ is not separable over $F$, for in that case every root of
$\Ph_n(x,t)$ is vacuously a root of $\partial_x \Ph_n(x,t)$ and is
therefore a multiple root.  Next, observe that $K[\xi_n, t] =
K[\xi_n]$ since $t=\fn{n}(\xi_n)$. Since $K[\xi_n, t]=K[\xi_n]$ is
finite as a $K[t]$-module, it cannot be a field; so $K[\xi]$ is
isomorphic to $K[x]$. Since $K[x]$ is normal, the same holds for
$K[\xi]$.  
\end{proof}

\subsection{Global specializations}\label{subsec:specialization}

Here we wish to clarify the nature of the specialization maps $\ee_n
\rightarrow \overline K$ associated with specializing $t$ to $t_0 \in
K$ as well as the relationship between the iterated
monodromy group $\M_\ph$ and its specializations $\M_{\ph,t_0}$.
We do so by defining a notion of
\emph{global specialization}.  Let $\oo_{\ee_\ph}$ be the
integral closure of $K[t]$ in $\ee_\ph$.  By integrality (and the
going up theorem), the maximal ideal $(t-t_0)$ of $K[t]$ extends to a
prime ideal $\m$ of $\oo_{\ee_\ph}$ such that $\m \cap K[t] =
(t-t_0)$.  The ring $\oo_{\ee_\ph}/ \m$ is integral over $K$, so is
actually a field.  Thus $\m$ is maximal, and $\oo_{\ee_\ph}/ \m$ is
algebraic over $K$.  So there are embeddings $\oo_{\ee_\ph} / \m
\rightarrow \overline K$. Fix one, and consider the associated map
$\sig:\oo_{\ee_\ph} \rightarrow \overline K$ with kernel $\m$. We call
such a map a \emph{global specialization} associated with $t_0$.  The
image of the global specialization, which is a field $\kk_{\ph,t_0}$, is
independent of the choice of global specialization $\sig$.

Now we can define the specializations $\sig_{n,t_0}:\oo_{\ee_n} \rightarrow
\overline K$ and $\oo_{F_n} \rightarrow \overline K$ by restriction of
the global specialization.  The field $\kk_{n, t_0}$ can be defined as
the image of $\oo_{\ee_n} \rightarrow \overline K$, and can be shown to
be independent of the choice of global specialization (associated with
$t_0$). However, $K_{n, t_0}$, the image of $\oo_{F_n} \rightarrow
\overline K$, depends on the global specialization $\sig$ as well as
on the choice of $\xi_n$.



In this optic, the relationship between the groups $\M_{\ph}=
\mathrm{Gal} (\ee_{\ph}/F)$ and the group $\M_{\ph, t_0} =
\mathrm{Gal} (\kk_{\ph}/K)$ is elucidated as follows.  Let $D_\m$ be
the \emph{decomposition group} associated to~$\m$ (consisting of the
elements of $\mathrm{Gal} (\ee_{\ph}/F)$ fixing the chosen maximal
ideal~$\m$ of $\oo_{\ee_\ph}$).  Then $D_\m$ acts on $\oo_{\ee_\ph} /
\m$, and therefore on $\kk_{\ph,t_0}$. Thus we get a homomorphism
$D_\m \rightarrow \mathrm{Gal} (\kk_{\ph,t_0}/K)$.  As usual, this is
a surjection, and if $t_0$ is not in the post-critical set, then it is
actually an isomorphism.  Thus, for $t_0 \in K \setminus \cP_\ph$,
$\M_{\ph, t_0}$ is isomorphic to a subgroup $D_\m$ of $\M_{\ph}$,
hence it too has an action on the rooted tree $T_\ph$.

\subsection{Dynamical systems on $\PP^1$}\label{subsec:dynamical}

\begin{definition}
Two self-maps $\ph,\psi$ of $\PP^1$ defined over $K$ (i.e. $\ph, \psi
\in K(x)$), are equivalent over $K$ (or $K$-conjugate) if there exists
an automorphism $\gamma$ of $\PP^1$ (defined over $K$) such that the
diagram
\begin{equation}
\begin{diagram}
\divide\dgARROWLENGTH by2
\node{\PP^1}
        \arrow[2]{e,t}{\gamma}
        \arrow{s,l}{\ph}
\node[2]{\PP^1}
        \arrow{s,r}{\psi}
\\
\node{\PP^1}
        \arrow[2]{e,t}{\gamma}
\node[2]{\PP^1}
\end{diagram}
        \label{fgisom}
\end{equation}
commutes.  In other words, $\ph$ and $\psi$ are equivalent over $K$ if
and only if there exist $a,b,c,d \in K$ satisfying $ad-bc\neq 0$ such
that $\ph(x)=\gamma^{-1}\psi\gamma(x)$ where $\gamma(x)=
\frac{ax+b}{cx+d}$.  The equivalence (or conjugacy) class of $\ph$,
denoted $[\ph]$, is a {\em dynamical system} on $\PP^1$.  For $\ph \in
\C(x)$, we say $[\ph]$ is {\em arithmetic} if there exists $\psi \in
\overline{\Q}(x)$ with $[\ph]=[\psi]$.
\end{definition}
Note that if $\ph\in K[x]$ is a polynomial map, the images of $\ph$
under {\em affine} transformations $\gamma(x)=ax+b$ over $K$ form
exactly the set of polynomial maps $K$-isomorphic to $\ph$.  Also, if
$\gamma$ takes $\fn{n}$ to $\psi^{\circ n}$ for $n=1$, then it does so
for all $n\geq 1$.  Thus, the study of iterations of $\ph$ and $\psi$
coincide (they simply take place in different coordinates) precisely
when $\ph$ and $\psi$ are conjugate.  In particular, if
$[\ph]=[\psi]$, then the iterated extensions $\ee_\ph$ and $\ee_\psi$
are isomorphic.  For a more detailed discussion, including the
relationship between fields of moduli and fields of definition of
dynamical systems on $\PP^1$, we refer the reader to Silverman
\cite{silverman}.

When discussing the coefficients of a post-critically finite
polynomial, it is often convenient to normalize by working with
\emph{monic} post-critically finite polynomials.

\begin{lemma} \label{lemma1}
Every polynomial $\ph\in K[x]$ of degree $d>1$ is equivalent over
some finite extension $K'/K$ to
a monic polynomial in $K'[x]$. Furthermore, if $\psi$ and $\ph$
are two $K'$-equivalent monic polynomials for some finite extension
$K'/K$,  then
$$
\psi(x) = \zeta^{-1} \ph(\zeta x + c) - \zeta^{-1} c
$$ where $c$ is in $K'$ and $\zeta$ is a $(d-1)$st root of unity.
\end{lemma}

\begin{proof}
Suppose $a x^d$ is the leading term of $\ph$.  If $\gamma(x) = b x +
c$, then $\gamma^{-1}(x) = b^{-1} x - b^{-1} c$. So $\gamma^{-1} \circ
\ph \circ \gamma(x)$ has leading term $b^{d-1} a x^d$. When we let $b$
be a root of $x^{d-1} - a^{-1}$, we find that $\gamma^{-1} \circ \ph
\circ \gamma$ is monic.  Now let $\ph$ and $\psi$ be monic equivalent
polynomials in $K'[x]$.  If $\gamma(x) = b x + c$, then $\gamma^{-1}
\circ \ph \circ \gamma(x)$ has leading term $b^{d-1} x^d$. Thus if
$\psi = \gamma^{-1} \circ \ph \circ \gamma$, then $b$ must be a
$(d-1)$th root of unity.
\end{proof}

\subsection{Examples: critically fixed simply ramified polynomials}
\label{subsec:cfsr}

Post-critically finite polynomials can be classified in terms of
certain combinatorial objects called {\em Hubbard trees}, see
\cite{bfh} and \cite{poirier}, as well as \cite{pilgrim} for their
relationship, in the case of two critical values, to {\em dessins
d'enfant} of genus 0.  Instead of describing this classification, in
this subsection, we simply want to illustrate that post-critically
finite polynomials are in plentiful supply by describing some of the
most simplest families of examples.
In order to avoid rationality questions, in this subsection we assume
that $K=\overline{K}$ is algebraically closed.  To write down
examples, we can make various simplifying assumptions; for example we
can limit the number of critical points (or values).  If $\ph$ has
only one critical point and this point is fixed, we see quickly that
$\ph$ is conjugate to $x\mapsto x^d$; specializations of this map
constitute the classical theory of ``pure'' extensions.  Another
family of examples is given by the {\em Chebyshev polynomials} which
have only two critical values; we study the quadratic one $x^2-2$ in
\S \ref{sec:quad}.  More generally, polynomials with two critical
values are called {\em generalized Chebyshev polynomials} or more
commonly {\em Shabat polynomials}; they have quite a rich structure,
as can be seen from the readable survey of Shabat-Zvonkin \cite{SZ}.

Here we make a different set of simplifying assumptions, and
completely classify the resulting post-critically finite dynamical
systems for each degree $d>1$.  Namely, we assume that the critical
points are fixed and that all the ramification indices are two; the
latter condition is equivalent to requiring that the polynomial has
$d-1$ critical points.  Other than $\ph(x)=x^d$, this is the simplest
family of post-critically finite polynomials.  It gives simple examples
of post-critically finite polynomials not equivalent to any monic
polynomial with {\em integer} coefficients.

\begin{definition}
A polynomial $\ph\in K[x]$ of degree~$d>1$ is said to be 
\emph{critically fixed, simply ramified} (CFSR) if 
$\ph$ has $d-1$ critical points, each of which is a fixed point
for $\ph$.
\end{definition}

\begin{example}
If $K$ does not have characteristic $2$, the polynomial $\ph(x)=x^2$
has exactly one critical point, $x=0$, which is a fixed point. Thus
$\ph$ is a CFSR polynomial.  It is easy to see that $\ph$ is the
unique such polynomial, up to equivalence, of degree~$2$.
\end{example}

\begin{example}
Let $K=\overline{\Q}$.  The polynomial $\ph(x)=x^3+\frac{3}{2}x$ has
derivative $\ph'(x)= 3 x^2 + \frac{3}{2}$.  Thus $\ph$ has two critical
points $\pm \frac{i}{\sqrt{2}}$.  The fixed points of $\ph$ are $0$ and
the two critical points, so $\ph$ is a CFSR polynomial.

This polynomial $\ph$ gives an example of a monic, post-critically
finite polynomial which does not have integral coefficients. Is there a
monic polynomial $\psi$ equivalent to $\ph$ with integer coefficients?  By
Lemma~\ref{lemma1} we only need to consider polynomials of the form
$$
\psi(x) = \ph(x + c) - c 
\quad\hbox{or}\quad
\psi(x) = - \ph(- x + c) + c.
$$ In the first case,
$$\psi(x) = (x+c)^3+\frac{3}{2}(x+c) -c
= x^3 + 3 c x^2 + \left( 3 c^2 + 
\frac{3}{2}\right) x + \left( c^3  + \frac{1}{2} c \right).
$$
Let $v$ be a place (valuation) in $\Q(c)$ above $2$
normalized so that $v(2)=1$.
We want to find $c$ so that the coefficients are integral.
So, $v\left( 3 c^2 + \frac{3}{2}\right) \ge 0$.
This implies $v(c) = -1/2$. 
Thus the coefficient of $x^2$ is not $2$-integral. A similar argument
applies to the second case.  We conclude that 
there are no monic polynomials with integral coefficients equivalent to~$\ph$.

This gives an example of a post-critically finite polynomial not equivalent
to any monic polynomial with integral coefficients.
\end{example}

We will now assume that $K$ is algebraically closed.  Thus, up to
equivalence, CFSR polynomials can be taken to be monic.  In an effort
to normalize further, consider the roots of the fixed point polynomial
$\ph(x) - x$.  These include all $d-1$ critical points (roots of
$\ph'$), but the polynomial is of degree~$d$. Thus there is a $d$th
root $r$; here, we allow $r$ to be one of the $d-1$ critical points if
$\ph(x)-x$ has a double root.  After conjugating by a translation
$\gamma$, we can assume that $r=0$.  In particular, $\ph(x) - x = d^{-1} x
\ph'(x)$.  This motivates the following.

\begin{definition}
A \emph{normalized CFSR polynomial} $\ph\in F[x]$ is 
a monic polynomial with $\ph(x) - x = d^{-1} x \ph'(x)$.
\end{definition}

The above argument gives the following.

\begin{lemma}
If $K$ is algebraically closed, then every CFSR polynomial 
is equivalent to a normalized CFSR polynomial.
\end{lemma}

Next, we will show that over an algebraically closed field,
there is, up to equivalence, a unique CFSR polynomial of
each degree. 

Assume $\ph\in F[x]$ is a normalized
CFSR polynomial of degree $d$. 
We rewrite $\ph(x) - x = d^{-1} x \ph'(x)$ as 
$$\ph(x) =  x+ d^{-1} x \ph'(x).$$
By differentiating this equation we get  $\ph'(x)  =  
1 + d^{-1}  \ph'(x) + d^{-1} x  \ph''(x)$,
so 
$$
\ph'  = \frac{d +  x  \ph''(x)}{d-1}
\quad\hbox{and}\quad
\ph =  x+ d^{-1} x \left( \frac{d +  x  \ph''(x)}{d-1} \right)
= \frac{d}{d-1} x + \frac{1}{d(d-1)} x^2  \ph''(x).
$$
Differentiating the first of these gives
$\ph''(x)  = \frac{1}{d-1}(\ph''(x) + x  \ph'''(x))$.
So if $d>2$, $\ph''(x)  = \frac{1}{d-2} x  \ph'''(x)$.
Thus
$$
\ph(x) = \frac{d}{d-1} x + \frac{1}{d(d-1)} x^2  \left(\frac{1}{d-2} x  \ph'''(x)  \right)
=\frac{d}{d-1} x + \frac{(d-3)!}{d!} x^3  \ph'''(x).
$$
Continuing in this manner, we get that the $n$th 
derivative $\ph^{(n)}(x)$ is $ \frac{1}{d-n} x \ph^{(n+1)}(x)$
if $n\le d$.
So, for $n\le d$.
$$
\ph =\frac{d}{d-1} x + \frac{(d-n)!}{d!} x^{n}  \ph^{(n)}(x).
$$
In particular, if $n=d$ then 
$$
\ph(x) =\frac{d}{d-1} x + \frac{1}{d!} x^{d}  \ph^{(d)}(x)
= \frac{d}{d-1} x +  x^{d}.
$$
This gives uniqueness.  Existence follows from
the fact that $\ph(x) = \frac{d}{d-1} x +  x^{d}$ satisfies
the equation 
$\ph(x) - x  = d^{-1} x \ph'(x)$ so is a normalized CFSR polynomial.

\begin{proposition}\label{prop:cfsr}
The polynomial $\frac{d}{d-1} x +  x^{d}$ is the unique
normalized CFSR polynomial of degree~$d$.
\end{proposition}

\begin{question}
Is it true that all post-critically finite polynomials over the
complex numbers are equivalent to a monic polynomial with algebraic
coefficients?
\end{question}
The answer to this question is known to be positive not just for
CFSR polynomials (by Proposition \ref{prop:cfsr}) but for all
critically fixed polynomials, by a theorem of Tischler \cite{tischler};
see Pakovich \cite{pakovich} for more on critically fixed polynomials.


\section{Characterization of finitely ramified iterated towers}\label{sec:pcf}

In this section, we prove Theorem \ref{thm:1}.  The main tool is a
polynomial version of the Riemann-Hurwitz genus formula.  Recall that
the resultant of two polynomials $P$ and $R$ in $K[x]$ satisfies
$$\Res_x(P,R) = a^{\deg_x(P)} \prod_{R(\theta)=0} P(\theta),$$
where the product is over the roots $\theta \in \overline{K}$ of $R$
(counted with multiplicity) and $a$ is the coefficient of the highest
power of $x$ appearing in $R$.

\begin{lemma}(Simon)\label{lem:simon}
Suppose $R(x,y)=A(y)x+B(y)$ where $A(y),B(y)\in K[y]$ are polynomials with
resultant $\Res_y(A(y),B(y))=\pm 1$.  If $Q(y)=\Res_x(P(x),R(x,y))$, then
$$\disc_y Q(y) = (\disc_x P(x))^{\deg_y R(x,y)} \Res_x(P(x),\disc_y R(x,y)).
$$
\end{lemma}
\begin{Proof}
A proof can be found in the thesis of D.~Simon \cite{simon}, see Proposition 
IV.2.2 as well as Remarque on p. 92.
\end{Proof}

Now we take a degree $d$ polynomial $\ph(x)= a_d x^d + a_{d-1}x^{d-1}
+ \cdots + a_1 x + a_0 \in K[x]$ with leading coefficient $a_d$, and
put $$D_n := \disc_x({\Ph_n(x,t)})\in F,$$ where $\Ph_n$ is defined by
(\ref{pn}).

\begin{proposition}\label{prop:recur}
Suppose $\ph$ is a degree $d$ polynomial in $K[x]$ with leading
coefficient $a_d$, and put $A=(-1)^{d(d-1)/2} d^d a_d^{d-1}$.  Then,
for $n\geq 1$, the discriminant $D_n=\disc_x (\fn{n}(x)-t)\in K[t]$
satisfies the recurrences
\begin{eqnarray*}
 D_{n+1} &=& {A}^{d^n} D_n^d \prod_{\beta \in \B_\ph}
\Ph_n(\beta,t)^{M_\beta}\\
&=& {A}^{d^n} D_n^d \prod_{r \in \RR_\ph}
\Ph_{n+1}(r,t)^{m_r},
\end{eqnarray*} 
where $m_r, M_\beta$ are the multiplicities defined by (\ref{mr}) and
(\ref{Mbeta}).
\end{proposition}
\begin{Proof}  We apply Lemma \ref{lem:simon} 
with $P(x)=\Ph_n(x,t)$, $R(x,y)= x- \ph(y)$; note that 
$Res_y(A(y),B(y))=\Res_y(1,-\ph(y))=1$.  We have \begin{eqnarray*}
Q(y) &=& \Res_x(\Ph_n(x,t),x-\ph(y)) \\
&=& \Ph_n(\ph(y),t) \\
&=& \Ph_{n+1}(y,t).
\end{eqnarray*}
Thus, by Lemma \ref{lem:simon}, \begin{eqnarray*}D_{n+1}&=&\disc_y Q(y)\\
&=&(\disc_x \Ph_n(x,t))^d
\Res_x(\Ph_n(x,t),\disc_y(x-\ph(y)))\\
&=& 
D_n^d A^{\deg_x \Ph_n(x,t)} \prod_{\beta\in \B_\ph} 
\Ph_n(\beta,t)^{M_\beta},
\end{eqnarray*}
where the factor $A^{\deg_x \Ph_n(x,t)}$ occurs because $A$ is the
leading coefficient of $\disc_y(x-\ph(y))$ when considered as a
polynomial in $x$.  We have $$\prod_{\beta\in \B_\ph} (\fn{n}(\beta)
-t)^{M_\beta} = \prod_{r\in \RR_\ph} (\fn{n}(\ph(r)) -t)^{m_r} = \prod_{r \in
\RR_\ph} (\fn{n+1}(r)-t)^{m_r}.$$  This completes the proof.
\end{Proof}
We now have the tools for proving Theorem \ref{thm:1}, but we first 
make a convenient definition.
\begin{definition}\label{def:S}
For a post-critically finite polynomial $\ph \in \oo_K[x]$ of degree
$d$ with leading coefficient $a_d$, and $t_0 \in K$, we let
$S_{\ph,t_0}$ be the set of real infinite places of $K$ together with
those finite ones which do not vanish on at least one of the
following: $d$, $a$, and $t_0-\nu$ as $\nu$ runs over ${\mathcal
P}_\ph$.
\end{definition}
\begin{proof}[Proof of Theorem \ref{thm:1}]
After a simple change of variables, we may assume that 
$\ph \in \oo_K[x]$.
By Lemma \ref{lem:integral} and Proposition \ref{prop:recur}, the
primes that ramify in $\ee_\ph/{\ee_1}$ are precisely those
corresponding to the postcritical set ${\mathcal P}_\ph$, i.e.
$\{(t-\nu) : \nu\in {\mathcal P}_\ph\}$.  Thus ${\ee_\ph}/F$ is finitely
ramified if and only if $\ph$ is post-critically finite.  Now suppose
$\ph$ is post-critically finite and fix $t_0 \in K$.  By Proposition
\ref{prop:recur}, for every $n\geq 1$, a place of $K_1$ which does not
lie over $S_{\ph,t_0}$ is unramified in $K_{n,t_0}/K_1$.  Hence
$K_{\ph,t_0}/K$ is finitely ramified and consequently, so is
${\kk_{\ph,t_0}}/K$.
\end{proof}

\begin{corollary}\label{cor:imr}
Let $K$ be a number field, $\ph \in \oo_K[x]$ a post-critically finite
polynomial of degree $d>1$.  For $t_0 \in K\setminus \cP_\ph$, 
the action of $\mathrm{Gal}(\overline{K}/K)$ on $\kk_{\ph,t_0}$
induces an iterated monodromy representation $\rho_{\ph,t_0}: G_{K,S}
\twoheadrightarrow \M_{\ph,t_0}$, where $S=S_{\ph,t_0}$.
\end{corollary}

\section{Polynomials with good reduction}\label{sec:wild}

Our aim here is to prove Theorem \ref{thm:wild1}.
Throughout this section, we suppose $K$ is a characteristic 0 field
equipped with an ultrametric valuation $v$ having valuation ring
$\oo_v=\{\alpha\in K : v(\alpha)\geq 0\}$; we assume that
$v(K^\times)=\Z$.  The residue field of $K$ with respect to $v$,
i.e. the reduction of $\oo_v$ modulo its maximal ideal $\m_v=
\{\alpha\in K : v(\alpha)> 0\}$, is denoted $k_v$.  We assume that
$k_v$ has positive characteristic $p>0$.

\begin{definition}
A polynomial $\ph= \sum_{j=0}^d a_j x^j\in K[x]$ has 
{\em good reduction} at $v$ if $$0=v(a_d)\leq v(a_j) \qquad \text{~for~}
1\leq j\leq d-1.$$  In other words, $\ph$ has good reduction when
it is $v$-integral with $v$-unital leading coefficient.
\end{definition}


\begin{lemma}  \label{lem:w1}
Suppose $K'$ is an algebraic extension of $K$ and fix an extension
$v'$ of $v$ to $K'$.  Let $\ph\in K[x]$ be a polynomial of degree $d\geq 2$
with good reduction at $v$.  If $\alpha \in K'$ has $v'(\alpha)<0$,
then $\alpha$ is not preperiodic for $\ph$.
\end{lemma}
\begin{Proof} 
Suppose $\beta\in K'$ has $v'(\beta)<0$.  Since the leading
coefficient of $\ph$ is a $v$-adic unit, there is a unique term in the
sum $\ph(\beta)=\sum_{j=0}^d a_j \beta^j$ with minimal valuation, namely
$a_d \beta^d$.  Since $v'$ is ultrametric, we have
$v'(\ph(\beta))=d\cdot v'(\beta)<v'(\beta)$.  Applying this principle to
$\alpha, \ph(\alpha), \fn{2}(\alpha), \ldots$, we obtain
$v'(\fn{n}(\alpha))=d^n v'(\alpha) \rightarrow -\infty$.  Thus, the
set $\{ \fn{n}(\alpha)\}$ cannot be finite since
$\{v'(\fn{n}(\alpha))\}$ is not finite.
\end{Proof}

\begin{lemma}\label{lemma:zero}
Suppose $\ph\in K[x]$ is a polynomial of degree $d$ divisible by $p$
(the residue characteristic of $v$), has good reduction at $v$, and is
post-critically finite.  Then the image of $\ph'$ in $k_v[x]$ is
identically $0$.
\end{lemma}
\begin{Proof}
Let $K'$ be a splitting field for $\ph'(x)$ over $K$, $v'$ an
extension of $v$ to $K'$ with valuation ring $\oo_{v'}=\{\alpha \in
K':v'(\alpha)\geq 0\}$ and residue field $k_{v'}$.  By Lemma
\ref{lem:w1}, $v'(r)\geq 0$ for every critical point $r\in \RR_\ph$.
We have $\ph'(x) = d \prod_{r\in\RR_\ph}(x-r)^{m_r}$.  Since the
critical points $r\in \RR_\ph$ are $v'$-integral, $\ph'(x)=p\psi(x)$
with $\psi(x)\in \oo_{v'}[x]$.  Thus, $\ph'(x)$ is identically zero in
$k_{v'}[x]$ and hence also in $k_v[x]$.
\end{Proof}

\begin{proof}[Proof of Theorem \ref{thm:wild1}]
For a polynomial $h \in \oo_v[x]$, let us write $\ord_p(h)=m$ if
$h(x)/p^m$ is in $\oo_v[x]$ but $h(x)/p^{m+1}$ is not.  By 
Lemma \ref{lemma:zero}, $\ord_p(\ph')\geq 1$.  
Using the product rule for differentiation and induction on $n$,
we have $\ord_p((\fn{n})') \geq 1 + \ord_p((\fn{n-1})')\geq n$.  We
have
$$
\disc_x(\Ph_n(x,t_0)) = \prod_{\theta\in \overline{K},\fn{n}(\theta)=t_0}
(\fn{n})'(\theta).
$$
Since $\deg(\fn{n})=d^n$, $p^{nd^n}$ divides $\disc_x(\Ph_n(x,t_0))$.  This
completes the proof.
\end{proof}

\begin{remark}
Via the example $\ph(x)=x^2-2$, we will see in \S \ref{sec:quad} that
the bound given Theorem \ref{thm:wild1} is sometimes met.  More
generally, suppose $K=\Q$, $\ph$ has degree $d=p$ and $t_0\in \Z$ is
such that for all $n\geq 0$, $K_{n+1,t_0}/K_{n,t_0}$ is Galois of
degree $p$.  We also assume that $p$ is totally ramified in
$K_{n,t_0}/\Q$ for all $n$.  These criteria are met, for example, for
the Chebyshev polynomial of degree $p$ and $t_0=0$ (giving the
cyclotomic $\Z_p$-extension of $\Q$).  Then $\ord_p (\disc~K_ n)
\leq n p^{n}$.  Here is how to see this.  In $K_{m+1}/K_m$, one knows
that $G_j$, the higher ramification group of (lower) index $j$, is
trivial once $j > p^{m+1}/ (p-1)$.  Hence,
$$\sum_j \# (G_j-1)< p^{m+1};$$ the left hand side of the above inequality
is equal to the valuation at the
prime above $p$ of the different of $K_{m+1}/K_m$.
One concludes by using the discriminant formula in a tower of extensions.
\end{remark}

\begin{example}  It is not difficult to write down
polynomials $\ph\in \Z[x]$ such that there is no wild ramification in
the iterated tower of $\ph$.  According to Theorem \ref{thm:wild1}, if
such a polynomial is monic, it will not be post-critically finite, so
the resulting iterated tower of function fields will be infinitely
ramified.  Here is a quadratic example.  Let $\ph(x)=x^2+x+\mu$ with
$\mu \in \Z$.  Then $\disc_x(\Ph_n(x,t))$ is {\em odd} for all $t\in \Z$
(for instance by Proposition \ref{prop:recur}).  However, $\ph$ is not
post-cricially finite.  Indeed, its only critical point is $r=-1/2$.
For $v=\ord_2$ the $2$-adic valuation of $\Z$, $v(r)=-1$ is negative,
hence by Lemma \ref{lem:w1}, $\ph$ is not post-critically finite.
\end{example}

\section{Prime decomposition in towers}\label{sec:decomposition}

In this section, $K$ is a number field.  We now describe, in terms of
certain graphs, how primes of $K$ not dividing the discriminant of
$\Ph_n(x,t_0)$ (assumed to be irreducible) decompose when we adjoin a
root of this polynomial.  A simple consequence of this description is
that no finite prime of $K$ splits completely in $K_{\ph,t_0}/K$.

We first set up some notation.  We assume $\ph \in \oo_K[x]$ is 
post-critically finite.  
Recall the notation from \S \ref{sec:intro} regarding $F_n=F(\xi_n)$.
Fixing $t_0 \in \oo_K$, let us assume that $\Ph_n(x,t_0)$ is
irreducible over $K$ for all $n\geq 1$ and choose a coherent system
$(\xi_n|_{t_0})$ of their roots so that $K_{n,t_0}=K(\xi_n|_{t_0})$.
For the rest of this section, we assume $\p$ is a prime of $\oo_K$
which is not in $S_{\ph,t_0}$ (see Definition \ref{def:S}).  For such
$\p$, the splitting of $\p$ in the ring of integers of $K_{n,t_0}$
coincides with the splitting of $\p$ in the ring
$\oo_K[\xi_n|_{t_0}]$; the latter factorization mirrors exactly the
factorization of the polynomial $\Ph_n(x,t_0)$ over the residue field
$\F_\p=\oo_K/\p$.

For example, the primes of degree 1 in $\oo_K[\xi_n|_{t_0}]$ which lie
over $\p$ correspond to the roots of $\fn{n}(x)-t_0$ over $\F_\p$,
i.e. the points in $\F_\p$ whose image under the $n$th iterate of
$\ph$ is the image $\overline{t}_0$ of $t_0$ in $\F_\p$.  A prime of
degree $k$ lying over $\p$ corresponds to a Galois orbit of $k$ points
defined over a degree $k$ extension of $\F_\p$ mapping to
$\overline{t}_0$ by $\fn{n}$.  Such data is conveniently summarized in
terms of certain directed graphs we now define.

For $k\geq 1$, let $\F_{\p,k}$ be a degree $k$ extension of the
residue field $\F_\p$.  We denote by $\Gamma_{\ph,\p,k}$ the following
directed graph: the vertices are the elements of $\F_{\p,k}$ and the
graph has a directed edge $v\rightarrow w$ if and only if $\ph(v)=w$.
After we choose an ordering $\lambda_1, \ldots, \lambda_{q}$ of the
elements of $\F_{\p,k}$, the adjacency matrix $A_{\ph,\p,k}$ of
$\Gamma_{\ph,\p,k}$ has $ij$ entry $1$ if $\ph(\lambda_i)\ =
\lambda_j$ and $0$ otherwise.  We write $\Gamma_{\ph,\p}$,
$A_{\ph,\p}$ for $\Gamma_{\ph,\p,1}$ and $A_{\ph,\p,1}$.

For calculations, it is useful to note that
$A_{\fn{n},\p,k}=A_{\ph,\p,k}^n$.  In other words, the in-degree of a
vertex $v$ in $\Gamma_{\fn{n},p,k}$ is the number of length $n$ paths on
$\Gamma_{\ph,\p,k}$ ending at $v$.  
For example, let $\overline{t}_0=t_0 +\p$ be the vertex corresponding
to the reduction of $t_0$ modulo $\p$.  Then the following quantities
all coincide:
\begin{itemize}
\item[(a)] the number of degree 1 primes of $\oo_K[\xi_n|_{t_0}]$ over $\p$,
\item[(b)] the in-degree of $\overline{t}_0$ on $\Gamma_{\fn{n},\p}$,
\item[(c)] the sum of the entries in the column of $A_{\ph,\p}^n$
corresponding to $\overline{t}_0$,
\item[(d)] the number of length $n$ paths on $\Gamma_{\ph,\p}$ ending
at $\overline{t}_0$.
\end{itemize} 

Note that, by (c) for example, there are at most $|\F_\p|=\mathbb{N}\p$
degree $1$ primes of $K_{n,t_0}$ lying over $\p$, hence $\p$ does not
split completely in $K_{\ph,t_0}/K$.

More generally, we can count the number of primes of any given degree
over $\p$ by taking into account the action of
$\Gal(\overline{\F}_{\p}/\F_\p)$.  Namely,
the graph $\Gamma_{\ph,\p,k}$ has the following additional structure:
each vertex is ``colored'', we will say {\em weighted}, by a positive
divisor $m$ of $k$ where $m$ is the exact degree of that vertex over
$\F_\p$.  Furthermore, every directed edge has the property that the
weight of the initial vertex is a multiple of the weight of the
terminal vertex. Also $\Gal(\F_{\p,k}/\F_\p)$ acts on the graph and
the weight of a vertex equals the size of its orbit under this action.

Summarizing the discussion, we have the following Proposition describing
prime decomposition in $K_{n,t_0}/K$ in terms of graphs.

\begin{proposition}\label{prop:splitting}
Suppose $\ph\in \o_K[x]$ is post-critically finite and that $t_0 \in
\oo_K$ is such that $\Ph_n(x,t_0)$ is irreducible over $K$ for all
$n\geq 1$.  Suppose $\p \subset \oo_K$ is not in $S_{\ph,t_0}$. Then,
for $k\geq 1$, the number of degree $k$ primes of $K_{n,t_0}$ lying
over $\p$ is $N/k$ where $N$ is the number of paths of length $k$ on
$\Gamma_{\ph,\p,k}$ which start with a vertex of weight $k$ and end at
$\overline{t}_0$, the weight 1 vertex corresponding to the image of
$t_0$ in $\F_\p$.
\end{proposition}

\begin{remark}
Alternatively, one could take the quotient graph of
$\Gamma_{\ph,\p,k}$ by identifying vertices which are in the same
orbit of $\Gal(\F_{\p,k}/\F_\p)$, and give a vertex in the new graph
the weight equal to the number of points identified.  Then the degree
$k$ primes of $K_{n,t_0}$ lying over $\p$ are in bijective
correspondence with the paths of length $n$ on the quotient graph
starting with a vertex of weight $k$ and ending at $\overline{t}_0$.
We should note that as long as $\p \not \in S_{\ph,t_0}$, the
decomposition of $\p$ in $K_{n,t_0}$ depends only on the residue of
$t_0$ modulo $\p$.
\end{remark}

For a fixed pair $(\p,k)$ and $n$ tending to infinity, each graph
$\Gamma_{\fn{n},\p,k}$ has $\mathbb{N}\p^k$ vertices and an equal number of
edges, hence is one of a finite number of graphs.  Therefore, the
sequence $\Gamma_{\fn{n},\p,k}$, $n=1,2,\ldots$ is always eventually
periodic.  In fact, it is relatively simple to describe exactly what
happens to the sequence of graphs in our situation.  Each connected
component of $\Gamma_{\ph,\p,k}$ consists of a unique cycle or
``loop'' with a number of ``arms'' emanating from it.  The minimal
period of the sequence $(\Gamma_{\fn{n},\p,k})$ is the lowest common
multiple of the length of the unique loop in each connected component
of $\Gamma_{\ph,\p,k}$ and the preperiod is the least common multiple
of the length of the longest arm in each connected component of
$\Gamma_{\ph,p,k}$.  All of these facts are easily verified and left
as amusing exercises for the reader.  A highly interesting question is
whether one can capture the graph-theoretical description of prime
decomposition in iterated extensions via appropriate zeta and
$L$-functions.  Here, we settle for a typical example as an
illustration.

\begin{example}\label{ex:tame}
Let $\ph(x)=x^2+i \in \Z[i]$.  Let $\p=(3+2i)$ be a prime of norm 13.
We map $\Z[i] \rightarrow \F_\p \simeq \F_{13}$ by sending $i \mapsto
8$, and list the elements of $\F_{13}$ as $0,1, 2,\ldots,12$.  We write
down the adjacency matrix $A_{\ph,\p}$ and draw the graph for $\ph$
and $\fn{2}$.  
$$
A_{\ph,\p} = \left(\begin{array}{ccccccccccccc}
 0& 0& 0& 0& 0& 0& 0& 0& 1& 0& 0& 0& 0\\
 
 0& 0& 0& 0& 0& 0& 0& 0& 0& 1& 0& 0& 0\\
 
 0& 0& 0& 0& 0& 0& 0& 0& 0& 0& 0& 0& 1\\
 
 0& 0& 0& 0& 1& 0& 0& 0& 0& 0& 0& 0& 0\\
 
 0& 0& 0& 0& 0& 0& 0& 0& 0& 0& 0& 1& 0\\
 
 0& 0& 0& 0& 0& 0& 0& 1& 0& 0& 0& 0& 0\\
 
 0& 0& 0& 0& 0& 1& 0& 0& 0& 0& 0& 0& 0\\
 
 0& 0& 0& 0& 0& 1& 0& 0& 0& 0& 0& 0& 0\\
 
 0& 0& 0& 0& 0& 0& 0& 1& 0& 0& 0& 0& 0\\
 
 0& 0& 0& 0& 0& 0& 0& 0& 0& 0& 0& 1& 0\\
 
 0& 0& 0& 0& 1& 0& 0& 0& 0& 0& 0& 0& 0\\
 
 0& 0& 0& 0& 0& 0& 0& 0& 0& 0& 0& 0& 1\\
 
 0& 0& 0& 0& 0& 0& 0& 0& 0& 1& 0& 0& 0
\end{array}
\right).
$$
 
The graph $\Gamma_{\ph,\p}$:

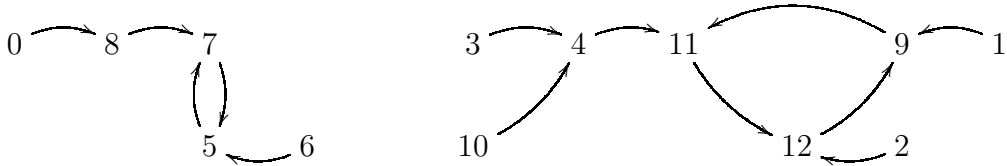
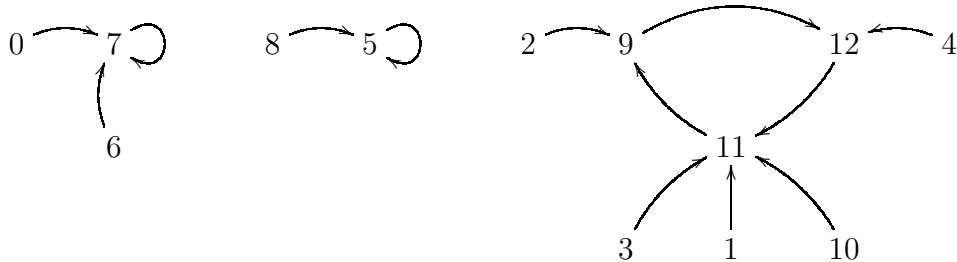
$$\xymatrix{ 0 \ar@/^/[r] & 8\ar@/^/[r] & 7 \ar@/^/[d] & \\ &&5
\ar@/^/[u]& \ar@/^/[l]6} \qquad \qquad \xymatrix{3 \ar@/^/[r]
&4\ar@/^/[r] &11 \ar@/_/[rd]&&\ar@/_0.5cm/[ll] 9&\ar@/_/[l] 1\\ 10
\ar@/_/[ur]&&&12 \ar@/_/[ur]&\ar@/^/[l]2& }$$

\

\

\

\

\

\

\

The graph $\Gamma_{\fn{2},\p}$:

$$\xymatrix{0\ar@/^/[r] &7\ar@(ur,dr)[]  \\
&6 \ar@/^/[u]} \qquad \qquad \xymatrix{8\ar@/^/[r] &5\ar@(ur,dr)[]} \qquad \qquad
\xymatrix{2\ar@/^/[r] &  9 \ar@/^0.5cm/[rr]&&12 \ar@/^/[ld]&\ar@/_/[l] 4\\
&&11 \ar@/^/[lu]&& \\
&3 \ar@/^/[ru]&1\ar[u]&10 \ar@/_/[lu] & }
$$



Note that $\Gamma_{\ph,\p}$ has two connected components, one with a
loop of length 2 and the other with a loop of length 3.  
The longest arm in each component has length 2.
The reader
can check either by taking powers of the adjacency matrix or by
drawing the graphs that $\Gamma_{\ph,\p}$ occurs only once in the
sequence $\Gamma_{\fn{n},\p}$, but starting with $n=2$, the sequence
has period $6$.  Note that $6$ is the product of the lengths of the
loops in the connected components of $\Gamma_{\ph,p}$.  With base
field $K=\F_\p$, the number of degree 1 places in $F_n$ over the prime
$(t-11)$ for $n=1,2,3,\ldots$ is the periodic sequence
$2,4,2,2,4,2,\ldots$ of period $3$.  As a check on the 
calculations, we verified using GP-PARI that with $\ph(x)=x^2+8$, 
the polynomials $\fn{n}(x)-11$ for $n=1,2,\ldots,7$, factor over $\F_{13}$ 
into distinct irreducible factors of the following degrees:
$$\begin{array}{ccc}
\hline
n & \text{~degrees of irreducible factors of~} \fn{n}(x)-11/\F_{13}& \text{no. of deg. 1 factors}\\
\hline
1& 1 ; 1&2\\ 
2& 1 ; 1 ; 1 ; 1&4 \\
3& 1 ; 1 ; 2 ; 2 ; 2&2 \\
4& 1 ; 1 ; 2 ; 2 ; 2 ; 2 ; 2 ; 4&2 \\
5& 1 ; 1 ; 1 ; 1 ; 2 ; 2 ; 2 ; 2 ; 4 ; 4 ; 4 ; 4 ; 4&4 \\
6& 1 ; 1 ; 2 ; 2 ; 2 ; 2 ; 2 ; 2 ; 2 ; 4 ; 4 ; 4 ; 4 ; 4 ; 4 ; 8 ; 8 ; 8&2 \\
7& 1 ; 1 ; 2 ; 2 ; 2 ; 2 ; 2 ; 2 ; 2 ; 2 ; 2 ; 4 ; 4 ; 4 ; 4 ; 4 ; 4 ; 4 ; 4 ; 4 ; 8 ; 8 ; 8 ; 8 ; 8 ; 8 ; 8 ; 8 ; 8&2 \\
\hline
\end{array}
$$
\end{example}

\section{Quadratic polynomials}\label{sec:quad}

In this section, 
we make a few remarks and give some examples concerning
quadratic polynomials.  By applying automorphisms of $\PP^1$, we bring
each quadratic polynomial to a standard form $\ph(x) = x^2 - r$.  We
then write down recurrence conditions for post-criticality of $\ph$.
The minimal number fields over which preperiodic points of prescribed
preperiod $m$ and period $n$ for such quadratic polynomials are
defined form an interesting family of number fields in their own
right.

\subsection{Normal form}
Put
$$\ph(x) = ax^2+bx+c \in K[x].$$
Let $\delta_\ph = -b^2/(4a) + c$.  It is the unique
branch point for the cover of $\PP^1$ given by the polynomial
$\ph(x)$, i.e. $\B_\ph = \{ \delta_\ph\}$.  
Theorem \ref{thm:1} now simplifies as follows: $\ee_\ph/F$
is finitely ramified if and only if $\delta_\ph$ is preperiodic for $\ph$.

If $\psi(x) = ax^2 + bx + c$ is quadratic, we
take $\gamma(x)=x/a$, so that $\gamma^{-1}(x)=ax$.  We then have
$$\gamma^{-1}\psi\gamma(x) = x^2 + bx + ac$$ is {\em monic}.  Note that
$\gamma$ fixes $0$.  Since an isomorphism from $\ph$ to $\psi$ carries
$\B_\ph$ to $\B_\psi$, applying a $K$-automorphism taking $\delta_\ph$ to
$0$, we see that $\psi$ is conjugate
to $\ph$ where
$$
\ph(x) = (x+\frac{b}{2})^2.
$$
We leave to the reader the exercise that for each quadratic $\psi\in
K[x]$, there is a unique $r\in \overline{K}$ such that $\psi$ is
conjugate to $(x-r)^2$.  Note that via the automorphism $\gamma(x)
=x+r$, the maps $x^2-r$ and $(x-r)^2$ are $K$-isomorphic.

Now consider a normalized quadratic polynomial $\ph(x)=x^2-r$.  We
have
$$\fn0(0)=0, \qquad \fn1(0) = -r , \qquad \fn2(0)=r (r-1), \qquad
\fn3(0)=r(r^3-2r^2+r-1), \cdots.$$ For $n\geq 0$, consider the
recurrence $g_{n+1}=rg_n^2-1$ with initial condition $g_0=0$.  Then
$\ph$ is post-critically finite if and only if $r$ is a root of
$g_m-g_n$ for some $m\neq n$.

\noindent{\bf Exercises.} i) If $r\in \Z$ and $\ph(x)=(x-r)^2$ has
periodic branch points, then $r \in \{0,1,2\}$.

ii) If $\ph(x)=ax^2+bx+c\in \overline{\Q}[x]$ has preperiodic branch
point, then $b/2$ is an algebraic integer.

\subsection{The polynomial $\ph(x)=x^2-2$}

In this subsection, we turn to an example which was the starting point
of this article.  While reading an article of Lemmermeyer, we came across
the classical fact that the cyclotomic $\Z_2$-extension of $\Q$ can be
written as $\Q(\theta_n)$ where $\theta_n= \sqrt{ 2+\sqrt{2+\cdots
+\sqrt{2}}})$.  Indeed, using the half-angle formula for cosines, one
establishes easily that the nested square root expression given above
evaluates $2\cos(\pi/2^{n+1})$.  What attracted our attention here was
that in the resulting recurrence-tower, the number of ramified primes
is finite (indeed only $2$ ramifies, and it does so totally and
deeply).  Since the $\theta_n$ are roots of the $n$th iterate equation
$\Ph_n(x,0)=\fn{n}(x)-0$, where $\ph(x)=x^2-2$, it was natural to
wonder whether for every $t\in K=\Q$, $\fn{n}(x)-t=0$ cuts out a
finitely ramified tower.  That this is so is guaranteed by Theorem
\ref{thm:1} since $x^2-2$ is post-critically finite.  Indeed, it is
the first member of the Chebyshev family of post-critically finite
polynomials.  For more details see Proposition 5.6 of \cite{bgn} where
the iterated monodromy group of any Chebyshev polynomial of degree
$d>1$ is shown to be infinite dihedral.

For the rest of this subsection, let $\ph(x)=x^2-2$.  Here we will
verify that another property of the cyclotomic $\Z_2$-tower (the
specialization of the tower at $t=0$) holds for many values of
$t_0\in\Z$, namely that a root of $\Ph_n(x,t_0)$ generates over $\Z$
the ring of integers of the number field it cuts out.

\begin{lemma}\label{lem:eisen}
For $t_0 \in \Z$, $t_0\equiv 0,1 \ \mod \ 4$, the polynomial
$\fn{n}(x)-t_0$ is irreducible over $\Q$.
\end{lemma}
\begin{proof}
We note that $\fn{n}(0)=-2$, $\fn{n}(\pm 1)=-1$.  If $t_0\equiv 0 \ \mod \ 4$,
we apply the Einsenstein criterion to $\fn{n}(x)$ at the prime $2$.  If
$t_0\equiv 1 \ \mod \ 4$, we use $\fn{n}(x+1)$ instead.
\end{proof}

\

\begin{proposition}
If $t_0 \in \Z$ is congruent to $0,1$ modulo $4$, and if $t_0+2$ and
$t_0-2$ are square-free, then for $n\geq 1$, the stem field
$K_n=\Q[x]/(\Ph_n(x,t_0))$ of the polynomial
$\Ph_n(x,t_0)=\fn{n}(x)-t_0$ is monog\`ene, as $\disc~K_n = \disc~
\Ph_n(x,t_0)$.
\end{proposition}
\begin{proof}
Letting $D_n= \disc( \Ph_n(x,t_0))$, by Proposition \ref{prop:recur},
we have for
$n\geq 1$,
$$ D_{n+1}= 4^{2^n} D_n^2 \Ph_n(-2),$$
or
\begin{equation}\label{recurD}
D_{n+1}= 4^{2^n} {D_n} ^2 (2-t_0)
\end{equation}
since $\fn{n}(\pm
2)=2$.  Also, for $n=1$, we have: $D_1= 4 (t_0+2)$.

We need to compare $D_n$ with the discriminant $d_n$ of the ring of
integers of $K_n$.  For $n=1$, we clearly have $d_n=D_n$, since
$t_0+2$ is square-free.  For $n\geq 1$, let us now determine the
ramification for each extension $K_{n+1}/K_n$.

We first remark that $K_{n+1}=K_n (\sqrt{\theta_n +2})$, with
$\Ph_n(\theta_n)=0$.  Next we observe that $N_{K_n/\Q} (\theta_n)=
\Ph_n(0)=\fn{n}(0)-t_0=2-t_0$.  Hence, for $n\geq 1$, in the extension
$K_{n+1}/K_n$, only the places dividing $2(2-t_0)$ are allowed to
ramify.  We first examine the tame ramification.  Suppose $l$ is a prime
divisor of $2-t_0$. Then $2+t_0 \equiv 4 \ \mod \ l$ and so $l$ is split
in $K_1/\Q$.
Let $l$ be an odd prime divisor of $t_0-2$. Since $N_{K_n/\Q}
(\theta_n)= 2-t_0$, there exists a prime $\mathfrak{l}_n$ of $K_n$ lying
over $l$ which is ramified in $K_{n+1}/K_n$.  In fact, there are two
primes over $l$ in $K_1$.  One of them is totally ramified in
$K_n/K_1$.  The other is unramified.  Therefore, the valuation
$v_{\mathfrak{l}_n}$ at the prime ideal $\mathfrak{l}_n$ of the
different of the extension $K_{n+1}/K_n$ is precisely $2-1=1$.

It remains to study the wild ramification.  For $n\geq 1$, let us put
$$\pi_n=\begin{cases} \theta_n & \text{~if~} t_0 \equiv 0 \pmod{4}\\
1+\theta_n & \text{~if~} t_0 \equiv 1 \pmod{4}.
\end{cases}
$$
We note that $2$
is ramified in $K_1/\Q$ and that $\pi_1$ 
is a uniformizer for the unique place $\p_1$ of $K_1$ lying over $2$.
We will proceed by induction.  Suppose, for some $n\geq 1$, that $2$
is totally ramified in $K_n/\Q$ and that 
$\pi_n$ is a uniformizer of the unique place
$\p_n$ of $K_n$ lying over $2$.  We claim that
$1+\pi_n$
is not a square
modulo ${\pi_n}^{2^{n+1}+1}$.
To see this, let
us suppose that $1+\pi_n$ is a square modulo $\pi_n^{2^{n+1}+1}$.  Since
the residue field is $\F_2$, we get, in the case $t_0 \equiv 1 \pmod{4}$,
$$1+\pi_n =(1+ a\pi_n)^2 \pmod{\pi_n^{2^{n+1}+1}},$$ with $a\in \Z_2$,
which is impossible.  Thus, for $t_0\equiv 1 \ \mod \ 4$, Kummer
theory tells us that $K_{n+1}/K_n$ is ramified at the unique place
above $2$.  For $t_0\equiv 0 \ \mod \ 4$, the argument is simpler,
since, in that case, the valution of $2+\pi_n$ at $\pi_n$ is the same as
that of $\theta_n$, namely $1$.  By Kummer theory, $K_{n+1}/K_n$ is
ramified at the unique place above $2$.  In conclusion, $K_{n+1}/\Q$
is totally ramified at $2$.

If $t_0\equiv \ 0 \ \mod \ 4$, it is clear that $\theta_{n+1}$ is a
uniformizer of the unique place of $K_{n+1}$ lying over $2$.  The
same holds for $1+\theta_{n+1}$ when $t_0\equiv 1 \ \mod \ 4$; note that
$N_{K_{n+1}/K_n} (1+\theta_{n+1})=-(\theta_n+1)$.  This completes
the induction step.

Next, let us calculate conductors.  Let $\sigma$ be a generator of the
Galois group $\Gal(K_{n+1}/K_n)$.  Then
$$({\sqrt{2+\theta_n}})^{\sigma-1} -1 =-2.$$ The valuation at $\p_{n+1}$
of $2$ is $2^{n+1}$.  Hence, the element $\sigma$ belongs to
$G_{2^{n+1}}$, but not to $G_{2^{n+1}+1}$ (we are using the higher
ramification groups in the lower numbering).  Consequently,
$$v_{\p_{n+1}} (\diff (K_{n+1}/K_n))=\sum_i (\# G_i-1)= 2^{n+1}$$
where $v_{\p_{n+1}}$ is the valuation at $\p_{n+1}$ and $\diff
(K_{n+1}/K_n)$ is the different of the extension $K_{n+1}/K_n$.

Now we are able to determine the discriminant of  $K_n/\Q$.
We have the recurrence formula
$$
\begin{array}{rcl}
\pm d_{n+1} & =& d_n^2 N_{K_{n+1}/\Q} \diff (K_{n+1}/K_n)\\
&=& d_n^2 (t_0-2) 2^{2^{n+1}}\\
&=& d_n ^2 (t_0-2) 4^{2^{n}},
\end{array}$$
which coincides up to sign with the recurrence (\ref{recurD}) for
$D_n$.  We also have the coincidence of initial conditions, $d_1=D_1$.
Since $D_n/d_n$ is a square, we conclude that $d_n=D_n$ for all $n$,
and so $\oo_{K_n}=\Z[\theta_n]$.
\end{proof}

\section{Iterated monodromy representations: questions}\label{sec:conjectures}

In this section, we discuss in a bit more detail conjectural and known
properties of iterated monodromy representations, especially as
compared with those of $p$-adic representations.  We also list a
number of open problems.

Let us first recall a conjecture of Fontaine and Mazur: If $K$ is a
number field and $S$ is a finite set of places of $K$ none of which
has residue characteristic $p$, then all finite-dimensional $p$-adic
representations of $G_{K,S}$ factor through a finite quotient (see
Conj.~5a of \cite{fm} as well as Kisin-Wortmann \cite{kw}).  On the
other hand, infinite tamely and finitely ramified extensions of number
fields do exist (and are in plentiful supply) thanks to the criterion
of Golod and Shafarevich, see e.g. Roquette \cite{roquette}.  Thus, at
least for certain pairs $K,S$, there is a sizeable portion of
$G_{K,S}$ which is predicted to be invisible to finite-dimensional
$p$-adic representations.

When $S$ contains all places above $p$, it is also expected, by a
conjecture of Boston \cite{boston-inv} (which we recall below), that
$p$-adic representations do not capture all of $G_{K,S}$.  Suppose
$\bar{\rho} : G_{K,S} \rightarrow \Gl_m(\F_p)$ is a residual
representation of $G_{K,S}$. By Mazur's theory of deformations,
there exists a universal ring $R(\rho)$ (local, noetherian
and complete) and a versal deformation $\rho : G_{K,S} \rightarrow
\Gl_m(R(\rho))$ such that $\bar{\rho}$ is the restriction of $\rho$.
Let $L=L_{\bar{\rho}}$ be the subfield of $K_S$ fixed by $\ker
\bar{\rho}$.  We put $H=H_{\bar{\rho}}=\Gal(M/L)$ where $M$ is the
maximal pro-$p$ extension of $L$ inside $K_S$.  If $S$ contains all
place above $p$ ($p$ odd, or for $p$ even we assume $K$ is totally
complex), then the cohomological dimension of $H$ is at most $2$.  The
purely group-theoretical Conjecture B of Boston \cite{boston-inv}
concerning the rank-growth of subgroups of $\Gl_m(R(\rho))$, then
implies the {\em non-injectivity conjecture} (\cite{boston-inv},
p. 91) to the effect that $\rho$ is never injective. In a certain
sense, one expects that $\rho$ forgets a non-trivial part of $H$.

How can one shed light on those sides of arithmetic fundamental groups
which are apparently not illuminated by the theory $p$-adic
representations?  As a counterpoint to the Fontaine-Mazur conjecture,
a conjecture of Boston \cite{boston-tr} asserts that infinite tame
quotients of $\G_{K,S}$ possess faithful actions on rooted trees.
Iterated monodromy groups are canonically equipped with such an action
\cite{nar}.  It is therefore natural to seek such representations via
specializations of iterated towers of post-critically finite
polynomials, in the wild case as well as in the tame case.  In the
wild case, it would be interesting to produce iterated monodromy
representations whose image does not have any infinite $p$-adic
analytic quotients.  Since very little is known about the structure of
infinite tamely and finitely ramified extensions of number fields, the
following question is of particular interest.

\begin{question}\label{ques:1}
Is there a number field $K$ and a rational function $\ph$ on $\P^1/K$
of degree $d>1$ as well as a specialization at $t_0\in K$ of
(\ref{pn}) such that
\begin{enumerate}
\item[i)] for each $n\geq 1$, $\Ph_n(x,t_0)$ is irreducible over $K$,
(i.e. $K_{n,t_0}=K(\xi_n|_{t_0}))$ is a field of degree $d^n$ over $K$),
\item[ii)] there is a {finite} set $S$ of places of $K$ such that $K_n/K$
is unramified outside $S$ for all $n\geq 1$, and such that 
\item[iii)] $S$ does not contain any primes dividing $d$?\end{enumerate}
\end{question}

By Theorem \ref{thm:1}, it is possible to fulfill ii) by taking $\ph$
to be a post-critically finite polynomial.  Satisfying i) is not too
difficult either, since we can arrange a place of $K$ to ramify
totally in $K_n$ (essentially an Eisenstein condition, see for example Lemma
\ref{lem:eisen}).  Condition iii) asks that $K_\ph/K$  be {\em tamely
ramified}.  It is not difficult to arrange i) and iii) simultaneously
by imposing congruence conditions (e.g. see Example \ref{ex:tame}).  
However, satisfying all conditions
together appears to be quite difficult.

A positive answer to Question \ref{ques:1} would provide, for the
first time, an explicit method for constructing an infinite tamely and
finitely ramified extension of a number field.  Currently the only
method for producing such extensions is via the Golod-Shafarevich
criterion.  On the other hand, a negative answer would assert that an
analogue of the Fontaine-Mazur Conjecture holds for finitely ramified
iterated extensions.  We should mention that for the function field of
a curve over a finite field with a square number of elements,
recursive constructions of Garcia-Stichenoth (see \cite{gar-st}
for example) for tamely
and finitely ramified extensions exist; that such constructions always arise
from modular curves is a conjecture of Elkies \cite{elkies}.

Recall that an algebraic extension $L$ over a number field $K$ is
called {\em asymptotically good} if i) $L/K$ is infinite, and ii) for
every sequence of distinct intermediate subfields of $L/K$, the root
discriminant\footnote{the root discriminant of a number field of
degree $n$ over $\Q$ is the $n$th root of the absolute value of its
discriminant} remains bounded.  A more general and more concise
version of Question \ref{ques:1} is the 
following.

\begin{question}\label{ques:2}
Is there a rational function $\ph$ on $\PP^1$ defined over a number
field $K$, and a $t_0 \in K$ such that the resulting specialized
iterated tower $K_{\ph,t_0}/K$ is asymptotically good?
\end{question}

Under the assumption of good reduction of the
polynomial $\ph$, the analogue  of this question where
we replace the number field discriminant with the polynomial
discriminant, has a negative answer by Theorem \ref{thm:wild1}.
Namely, for a polynomial $P \in \Q[x]$ of degree $d\geq 1$, define its root
discriminant by $\rd(P)=|\disc(P)|^{1/d}$.  An immediate consequence
of Theorem \ref{thm:wild1} is
\begin{corollary}\label{cor:wild1}
If $\ph\in \Q[x]$ is post-critically finite, has degree divisible by
$p$, and has good reduction at $p$, then for any $t_0 \in \Z$, the
sequence of polynomials $(\Ph_n(x,t_0))$ is asymptotically bad in the
sense that (the $p$-part of) $\mathrm{rd}(\Ph_n(x,t_0))$ tends to
infinity with $n$.
\end{corollary}
This result is in agreement with a conjecture of Simon \cite{simon},
to the effect that any infinite sequence of distinct polynomials over
$\Z$ is asymptotically bad.  Thus, to tackle Questions \ref{ques:1} and
\ref{ques:2}, one would very likely have to understand the index of 
the order $\oo_K[\xi_n|_{t_0}]$ in $\oo_{K_{n,t_0}}$.




\begin{thebibliography}{99}



\bibitem[BGN]{bgn} L.~Bartholdi, R.I.~Grigorchuk and V.~Nekrashevych, From
fractal groups to fractal sets, preprint 
\verb?arXiv:math.GR/0202001?, 75pp.

\bibitem[BORT]{bort} H. Bass, M. Otero-Espinar, D. Rockmore, and C. Tresser,
Cyclic renormalization and automorphism groups of rooted trees.
Lecture Notes in Mathematics, 1621.
Springer-Verlag, Berlin, 1996. 
MR1392694 (97k:58058) 


\bibitem[BFH]{bfh} B. Bielefeld, Y. Fisher, and J. Hubbard, 
The classification of critically preperiodic polynomials as dynamical systems.
J. Amer. Math. Soc. 5 (1992), no. 4, 721--762. MR1149891 (93h:58128)

\bibitem[B]{boston-inv} N.~Boston, Explicit deformation of Galois
representations.  Invent. Math.  103 (1991), no. 1, 181--196.
MR1079842 (91j:11041)

\bibitem[B1]{boston-tr} N.~Boston, Tree representations of Galois
groups, preprint, \verb?www.math.wisc.edu/~boston?, 5pp.



\bibitem[BP]{bux-perez} K-U.~Bux and R.~Perez, 
On the growth of iterated monodromy groups, preprint 
\verb?arXiv:math.GR/? \verb?0405456?, 15pp.


\bibitem[DH]{dh} A.~Douady and J.~Hubbard, A proof of Thurston's
topological characterization of rational functions.  Acta Math.  171
(1993), no. 2, 263--297 MR1251582 (94j:58143)


\bibitem[E]{elkies} N.~Elkies, Explicit modular towers, pages 23-32 in
Proceedings of the Thirty-Fifth Annual Allerton Conference on
Communication, Control and Computing (1997, T. Basar, A. Vardy, eds.),
Univ. of Illinois at Urbana-Champaign 1998 (\verb?arXiv:math.NT/0103107?)



\bibitem[FM]{fm} J.-M.~Fontaine and B. Mazur, Geometric Galois
 representations. Elliptic curves, modular forms, \& Fermat's last
 theorem (Hong Kong, 1993), 41--78, Ser. Number Theory, I,
 Internat. Press, Cambridge, MA, 1995

\bibitem[GSR]{gar-st} A. Garcia, H. Stichtenoth, and H-G. R\"uck, On tame
towers over finite fields. J. Reine
Angew. Math. 557 (2003), 53--80.  MR1978402 (2004e:11133)

\bibitem[KW]{kw} M.~Kisin and S.~Wortmann, A note on Artin motives.
Math. Res. Lett.  10 (2003), no. 2-3, 375--389.  MR1981910
(2004d:14018)


\bibitem[N]{nar} V. Nekrashevych, Iterated Monodromy Groups, preprint
\verb?arXiv:math.DS/0312306?, 40pp.

\bibitem[P]{pakovich} F. Pakovich, Arithmetics of conservative
    polynomials and yet another action of $\Gal(\bar \Q/\Q)$ on plane
    trees, preprint \verb?arXiv:math.NT/0404478?, 9pp.

\bibitem[P]{pilgrim} K. Pilgrim,  Dessins d'enfants and Hubbard
trees.  Ann. Sci. \'Ecole Norm. Sup. (4) 33 (2000), no. 5,
671--
693. MR1834499 (2002m:37062)

\bibitem[Po]{poirier} A. Poirier, On postcritically finite polynomials,
part 1 \verb?arXiv:math.DS/9305207? 45pp., and part 2
\verb?arXiv:math.DS/9307235? 36pp.

\bibitem[R]{roquette} P. Roquette, On class field towers. 
Algebraic Number Theory (Proc. Instructional Conf., Brighton, 1965)
pp. 231--249 Thompson, Washington, D.C. 1967. MR0218331 (36 \#1418)


\bibitem[Si]{silverman} J.H.~Silverman, The field of definition
for dynamical systems on $\mathbf P\sp 1$.  Compositio Math.  98 (1995),
no. 3, 269--304.  MR1351830 (96j:11090)


\bibitem[S]{simon} D.~Simon, \'Equations dans les corps de nombres et
discriminants minimaux, Th\`ese, Universit\'e de Bordeaux I, 21 Dec. 1998.
\verb?http://www.math.unicaen.fr/~simon/?

\bibitem[SZ]{SZ} G. Shabat and A. Zvonkin, Plane trees and algebraic
Numbers.  Jerusalem combinatorics '93, 233--275, Contemp. Math., 178,
Amer. Math. Soc., Providence, RI, 1994. MR1310587 (96d:14028)


\bibitem[T]{tischler} D. Tischler, Critical points and values of complex 
polynomials, J. of Complexity 5 (1989), 438--456.  MR1028906 (91a:30004)


\end{thebibliography}
\end{document}